\documentclass[reqno]{amsart}
\usepackage{amsmath,amsrefs,amssymb}

\numberwithin{equation}{section}

\DeclareMathOperator{\divergence}{div}

\DeclareMathOperator{\Vol}{Vol}
\DeclareMathOperator{\Scal}{S}
\DeclareMathOperator{\Ricci}{Ric}
\DeclareMathOperator{\Einstein}{E}
\DeclareMathOperator{\Q}{Q}
\DeclareMathOperator{\Paneitz}{P}

\newcommand{\N}{\mathbb{N}}

\renewcommand{\[}{\left[}
\renewcommand{\]}{\right]}
\renewcommand{\(}{\left(}
\renewcommand{\)}{\right)}

\newtheorem{theorem}{Theorem}[section]

\begin{document}

\title[Uniqueness of conformal metrics with constant $\Q$-curvature]{Uniqueness of conformal metrics  with constant $\Q$-curvature on closed Einstein manifolds}

\author{J\'er\^ome V\'etois}

\address{J\'er\^ome V\'etois, Department of Mathematics and Statistics, McGill University, 805 Sherbrooke Street West, Montreal, Quebec H3A 0B9, Canada}
\email{jerome.vetois@mcgill.ca}

\thanks{The author was supported by the NSERC Discovery Grant RGPIN-2022-04213.}

\date{October 13, 2022}

\thanks{The author would like to thank very much Robin Graham, Matthew Gursky, Emmanuel Hebey and Fr\'ed\'eric Robert for helpful advice and suggestions as well as Matthew Gursky and Andrea Malchiodi for very interesting comments on the manuscript.}

\begin{abstract}
On a smooth, closed Riemannian manifold $\(M,g\)$ of dimension $n\ge3$ with positive scalar curvature and not conformally diffeomorphic to the standard sphere, we prove that the only conformal metrics to $g$ with constant $\Q$-curvature of order 4 are the metrics $\lambda g$ with $\lambda>0$ constant.
\end{abstract}

\maketitle

\section{Introduction and main result}\label{Sec1}

On a smooth, closed (i.e. compact and without boundary) Riemannian manifold of dimension $n\ge3$, Branson's $\Q$-curvature~\cite{Bra1} is defined as
\begin{equation}\label{Sec1Eq1}
\Q_g:=\frac{1}{2\(n-1\)}\Delta_g\Scal_g+\frac{n^3-4n^2+16n-16}{8\(n-1\)^2\(n-2\)^2}\Scal_g^2-\frac{2}{\(n-2\)^2}\left|\Ricci_g\right|_g^2,
\end{equation}
where $\Delta_g:=-\divergence_g\nabla$ is the Laplace--Beltrami operator, $\Scal_g$ is the scalar curvature and $\Ricci_g$ is the Ricci curvature of the manifold. 

\medskip
In the case of the standard sphere, the conformal metrics with constant $\Q$-curvature have been classified by Lin~\cite{Lin} by using the moving-plane method. In this case, there exists an explicit, multi-dimensional family of conformal metrics with constant $\Q$-curvature. In this article, we examine the case of Einstein manifolds with positive scalar curvature and not conformally diffeomorphic to the standard sphere. In this case, we obtain the following:

\begin{theorem}\label{Th1}
Let $\(M,g\)$ be a smooth, closed Einstein manifold of dimension $n\ge3$ with positive scalar curvature and not conformally diffeomorphic to the standard sphere. Then the only conformal metrics to $g$ with constant $\Q$-curvature are the metrics $\lambda g$ with $\lambda>0$ constant. 
\end{theorem}

This result extends to the $\Q$-curvature problem a result obtained by Obata~\cite{Oba} for the scalar curvature problem (see also the subsequent articles by Bidaut-V\'eron and V\'eron~\cite{BidVer} and Gidas and Spruck~\cite{GidSpr} for extensions of Obata's result to more general second-order equations). 

\medskip
The transformation law for the $\Q$-curvature under a conformal change of metric is given by the equations
\begin{equation}\label{Sec1Eq2}
\left\{\begin{aligned}&\Paneitz_gu+\Q_g=\Q_{e^{2u}g}e^{4u}\quad\text{in }M&&\text{ if }n=4\\
&\Paneitz_gu=\frac{n-4}{2}\Q_{u^{\frac{4}{n-4}}g}u^{\frac{n+4}{n-4}},\ u>0\quad\text{in }M&&\text{ if }n\ne4,\end{aligned}\right.
\end{equation}
where $\Paneitz_g$ is the Paneitz--Branson operator~\cites{Bra1,Pan} defined as
$$\Paneitz_g:=\Delta_g^2-\divergence_g\(\(\frac{n^2-4n+8}{2\(n-1\)\(n-2\)}\Scal_gg-\frac{4}{n-2}\Ricci_g\)\nabla u\)+\frac{n-4}{2}\Q_g.$$
In particular, in the case where $\(M,g\)$ is Einstein, we obtain
$$Q_g=\frac{\(n+2\)\(n-2\)}{8n\(n-1\)^2}\Scal_g^2\quad\text{and}\quad \Paneitz_g=\Delta_g^2+\frac{n^2-2n-4}{2n\(n-1\)}\Scal_g\Delta_g+\frac{n-4}{2}\Q_g.$$
We mention in passing that the concept of $Q$-curvature and the corresponding operator $P_g$ have been shown to have natural extensions to higher orders. Some references in this case are by Branson~\cites{Bra2}, Fefferman and Graham~\cites{FefGra1,FefGra2}, Gover~\cite{Gov}, Graham, Jenne, Mason and Sparling~\cite{GraJenMasSpa} and Juhl~\cite{Juhl}. In particular, in the case of Einstein manifolds, explicit formulas can be found for all orders (see \cite{FefGra2}*{Proposition~7.9} and also \cite{Gov}).

\smallskip
We use a similar approach as in Obata's proof~\cite{Oba} for the scalar curvature problem (see also \cites{BidVer,GidSpr}). Given a conformal metric $g_v:=v^{-1}g$, where $v\in C^\infty\(M\)$, $v>0$ in $M$, this approach consists in finding a suitable function $\Theta^{\(k\)}_g\(v\)\in C^\infty\(M\)$, with $k=1$ for the scalar curvature problem and $k=2$ for the $\Q$-curvature problem, such that if $\Scal_{g_v}$ is constant for $k=1$ and $\Q_{g_v}$ is constant for $k=2$, then 
$$\Theta^{\(k\)}_g\(v\)\ge0\text{ in }M,\ [\Theta^{\(k\)}_g\(v\)\equiv0\text{ in }M\,\Longleftrightarrow\,\nabla v\equiv0\text{ in }M]\text{ and}\int_M\Theta^{\(k\)}_g\(v\)dv_g=0.$$
The last equality is obtained by applying multiple integrations by parts. It easily follows from these properties that if such a function $\Theta^{\(k\)}_g\(v\)$ exists, then $v$ must be constant. For $k=1$, this is achieved by considering the function
$$\Theta^{\(1\)}_g\(v\):=v^{\frac{3-n}{2}}\left|\Einstein_{g_v}\right|_g^2,$$
where $\Einstein_{g_v}$ is the Einstein tensor of the metric $g_v$, i.e.
\begin{equation}\label{Sec1Eq3}
\Einstein_{g_v}:=\Ricci_{g_v}-\frac{1}{n}\Scal_{g_v}g_v=\Ricci_g-\frac{1}{n}\Scal_gg+\(n-2\)\sqrt{v}^{-1}\Big(\nabla^2\sqrt{v}+\frac{1}{n}\Delta_g\sqrt{v}\,g\Big).
\end{equation}
For $k=2$, we use the function
\begin{multline}\label{Sec1Eq4}
\Theta^{\(2\)}_g\(v\):=v^{\frac{1-n}{2}}\Bigg(\left|\nabla\Scal_{g_v}+\frac{3n-4}{2\(n-2\)}\Einstein_{g_v}\nabla v\right|^2_g-\frac{\(3n-4\)^2}{4\(n-2\)^2}\big|\Einstein_{g_v}\nabla v\big|_g^2\\+\frac{\left|\Einstein_{g_v}\right|^2_g}{\(n-2\)^2}\big(2\(n^2-2\)\Scal_{g_v}v+4\(n-1\)\Scal_gv^2+n\(n-1\)^2\left|\nabla v\right|_g^2\big)\Bigg).
\end{multline}
Under the conditions of Theorem~\ref{Th1}, we can see that the function $\Theta^{\(2\)}_g\(v\)$ is non-negative in $M$ by observing that 
\begin{equation}\label{Sec1Eq5}
\(3n-4\)^2\big|\Einstein_{g_v}\nabla v\big|_g^2\le\(3n-4\)^2\left|\Einstein_{g_v}\right|_g^2\left|\nabla v\right|_g^2\le 4n\(n-1\)^2\left|\Einstein_{g_v}\right|_g^2\left|\nabla v\right|_g^2
\end{equation}
and using the positivity of the functions $v$, $\Scal_g$ and $\Scal_{g_v}$, the latter following from a result obtained by Gursky and Malchiodi~\cite{GurMal} (see Theorem~\ref{ThGM} below). Finally, it is not difficult to see that $\Theta^{\(2\)}_g\(v\)\equiv0$ in $M$ if and only if $\nabla\Scal_{g_v}\equiv0$ in $M$ and either $\Einstein_{g_v}\equiv0$ or $\nabla v\equiv0$ in $M$, which in both cases, is equivalent to $v$ being constant, provided $\(M,g\)$ is Einstein and not conformally equivalent to the standard sphere. In the case where $\(M,g\)$ is conformally equivalent to the standard sphere, this approach also gives an alternative proof of the classification of conformal metrics with constant $Q$-curvature.

\smallskip
In the case of more general manifolds, the problem \eqref{Sec1Eq2} has been studied my several authors. Existence results have been obtained by Brendle~\cite{Bre}, Chang and Yang~\cite{ChangYang}, Djadli and Malchiodi~\cite{DjaMal} and Li, Li and Liu~\cite{LiLiLiu} in dimension $n=4$, by Djadli, Hebey and Ledoux~\cite{DjaHebLed}, Esposito and Robert~\cite{EspRob}, Gursky, Hang and Lin~\cite{GurHangLin}, Gursky and Malchiodi~\cite{GurMal}, Hang and Yang~\cites{HangYang1,HangYang3} and Qing and Raske~\cite{QingRas2} in dimensions $n\ge5$ and by Hang and Yang~\cite{HangYang2} in dimension $n=3$. Non-uniqueness results have been obtained by Bettiol, Piccione and Sire~\cite{BetPicSire} in dimensions $n\ge5$. The question of compactness of the set of solutions has also been studied by Druet and Robert~\cite{DruRob}, Malchiodi~\cite{Mal} and Weinstein and Zhang~\cite{WeiZhang} in dimension $n=4$, by Hebey and Robert~\cite{HebRob}, Li~\cite{Li}, Li and Xiong~\cite{LiXio} and Qing and Raske~\cite{QingRas1} in dimensions $n\ge5$ and by Hang and Yang~\cite{HangYang2} in dimension $n=3$. Finally non-compactness results have been obtained in dimensions $n\ge25$ by Wei and Zhao~\cite{WeiZhao}.

\section{Extended versions and proof of Theorem~\ref{Th1}}\label{Sec2}

In this section, we prove the following results, which are slightly more general than Theorem~\ref{Th1}:

\begin{theorem}\label{Th2}
Let $p\le4$ and $\(M,g\)$ be a smooth, closed Einstein manifold of dimension $n=4$ with positive scalar curvature. In the case where $p=4$, assume that $\(M,g\)$ is not conformally diffeomorphic to the standard sphere. Then there does not exist any non-constant solutions to the equation 
\begin{equation}\label{Th2Eq1}
\Paneitz_gu+\Q_g=e^{pu}\quad\text{in }M.
\end{equation}
\end{theorem}

\begin{theorem}\label{Th3}
Let $n\ge3$, $n\ne4$, $p\le\frac{2n}{n-4}$ and $\(M,g\)$ be a smooth, closed Einstein manifold of dimension $n$ with positive scalar curvature. In the case where $p=\frac{2n}{n-4}$, assume that $\(M,g\)$ is not conformally diffeomorphic to the standard sphere. Then there does not exist any non-constant solutions to the equation 
\begin{equation}\label{Th3Eq1}
\Paneitz_gu=u^{p-1},\ u>0\quad\text{in }M.
\end{equation}
\end{theorem}

We obtain Theorem~\ref{Th1} from Theorems~\ref{Th2} and~\ref{Th3} as follows:

\proof[Proof of Theorem~\ref{Th1}]
Let $\tilde{g}$ be a conformal metric to $g$ with constant $\Q$-curvature. Let $u\in C^\infty\(M\)$ be such that $\tilde{g}=e^{2u}g$ in $M$ if $n=4$ and $\tilde{g}=u^{\frac{4}{n-4}}g$ with $u>0$ in $M$ if $n\ne4$. Since $\Q_g$ and $\Q_{\tilde{g}}$ are constant, integrating \eqref{Sec1Eq2} in $M$ gives
\begin{equation}\label{Th1Eq1}
\Q_g\Vol_g\(M\)=\left\{\begin{aligned}&\Q_{\tilde{g}}\int_Me^{4u}dv_g&&\text{if }n=4\\&\Q_{\tilde{g}}\int_Mu^{\frac{n+4}{n-4}}dv_g&&\text{if }n\ne4,\end{aligned}\right.
\end{equation}
where $\Vol_g\(M\)$ is the volume of $\(M,g\)$ and $dv_g$ is the volume element of $\(M,g\)$. Since $Q_g>0$ in $M$, it follows from \eqref{Th1Eq1} that $\Q_{\tilde{g}}>0$ in $M$. We then define
$$\tilde{u}:=\left\{\begin{aligned}&u+\frac{1}{4}\ln\Q_{\tilde{g}}&&\text{if }n=4\\&\(\frac{n-4}{2}\Q_{\tilde{g}}\)^{\frac{n-4}{8}}u&&\text{if }n\ne4.\end{aligned}\right.$$
We can then rewrite \eqref{Sec1Eq2} as
$$\left\{\begin{aligned}&\Paneitz_g\tilde{u}+\Q_g=e^{4\tilde{u}}\quad\text{in }M&&\text{ if }n=4\\
&\Paneitz_g\tilde{u}=\tilde{u}^{\frac{n+4}{n-4}},\ \tilde{u}>0\quad\text{in }M&&\text{ if }n\ne4,\end{aligned}\right.$$
Applying Theorems~\ref{Th2} and~\ref{Th3}, we then obtain that $\tilde{u}$ is constant, which implies that $u$ is constant. This ends the proof of Theorem~\ref{Th1}.
\endproof

\smallskip
The proofs of Theorems~\ref{Th2} and~\ref{Th3} use the following result, which is a straightforward variation of a result obtained by Gursky and Malchiodi~\cite{GurMal}:

\begin{theorem}\label{ThGM}
Let $\(M,g\)$ be a smooth, closed Riemannian manifold of dimension $n\ge3$ with positive scalar curvature and non-negative $\Q$-curvature. Let $\tilde{g}$ be a conformal metric to $g$ with non-negative $\Q$-curvature. Then the scalar curvature of $\tilde{g}$ is positive.
\end{theorem}

\proof[Proof of Theorem~\ref{ThGM}]
We closely follow the proof of Theorem~2.2 in~\cite{GurMal}. Let $u\in C^\infty\(M\)$ be such that $\tilde{g}=e^{2u}g$ in $M$ if $n=4$ and $\tilde{g}=u^{\frac{4}{n-4}}g$ with $u>0$ in $M$ if $n\ne4$. Since $\Scal_g>0$ in $M$, we can define
$$t_0:=\sup\left\{t\in\[0,1\]:\,\Scal_{g_s}>0\text{ for all }s\in\[0,t\]\right\},$$
where $g_s:=e^{2su}g$ if $n=4$ and $g_s=\(1-s+su\)^{\frac{4}{n-4}}g$ in $M$ if $n\ne4$. Notice that by continuity, we obtain $\Scal_{g_{t_0}}\ge0$ in $M$. On the other hand, since $\Q_g\ge0$ and $\Q_{\tilde{g}}\ge0$ in $M$, using \eqref{Sec1Eq2}, we obtain  
\begin{equation}\label{ThGMEq1}
\Q_{g_{t_0}}=\left\{\begin{aligned}&e^{-4t_0u}\(\(1-t_0\)\Q_g+t_0\Q_{\tilde{g}}e^{4u}\)&&\text{if }n=4\\
&\(1-t_0+t_0u\)^{-\frac{n+4}{n-4}}\(\(1-t_0\)\Q_g+t_0\Q_{\tilde{g}}u^{\frac{n+4}{n-4}}\)&&\text{if }n\ne4\end{aligned}\right\}\ge0\quad\text{in }M.
\end{equation}
It follows from \eqref{Sec1Eq1} and \eqref{ThGMEq1} that 
\begin{equation}\label{ThGMEq2}
\Delta_{g_{t_0}}\Scal_{g_{t_0}}\ge-\frac{n^3-4n^2+16n-16}{4\(n-1\)\(n-2\)^2}\Scal_{g_{t_0}}^2\quad\text{in }M.
\end{equation}
Since $\Scal_{g_{t_0}}\ge0$ in $M$, by the strong maximum principle, it follows from \eqref{ThGMEq2} that either $\Scal_{g_{t_0}}\equiv0$ or $\Scal_{g_{t_0}}>0$ in $M$. Since $\Scal_g>0$ in $M$ and a conformal class cannot contain both a metric with positive scalar curvature and a metric with zero scalar curvature, we then obtain that $\Scal_{g_{t_0}}>0$ in $M$, which implies that $t_0=1$ and $\Scal_{\tilde{g}}>0$ in $M$. This ends the proof of Theorem~\ref{ThGM}.
\endproof

Let us now set some notations and recall some preliminary formulas. We let $\(\cdot,\cdot\)_g$ be the multiple inner product induced by the metric $g$ for the tensors of same rank (i.e. such that for example $\(S,T\)_g=g^{i_1j_1}\dotsm g^{i_lj_l}S_{i_1\dotsc i_l}T_{j_1\dotsc j_l}$
for covariant tensors $S$ and $T$ of rank $l\in\N$). We let $\left|\cdot\right|_g$ be the norm induced by $\(\cdot,\cdot\)_g$. We denote by $\overline{\Delta}_g$ the connection Laplacian and by $\Delta_g$ the Hodge Laplacian on 1-forms. Given a smooth function $u$ in $M$, the Weitzenbock identity gives
\begin{equation}\label{Sec2Eq1}
\Delta_gdu=\overline{\Delta}_gdu+\Ricci_g\nabla u.
\end{equation}
We also recall the Bochner--Lichnerowicz--Weitzenbock formula 
\begin{equation}\label{Sec2Eq2}
\frac{1}{2}\Delta_g\left|\nabla u\right|_g^2=\(\nabla\Delta_g u,\nabla u\)_g-\left|\nabla^2 u\right|_g^2-\Ricci_g \(\nabla u,\nabla u\).
\end{equation}

\medskip
We begin with proving Theorem~\ref{Th2}.

\proof[Proof of Theorem~\ref{Th2}]
Let $u$ be a solution of \eqref{Th2Eq1}. Using \eqref{Sec2Eq1} together with the fact that $\(M,g\)$ is Einstein, we then obtain
\begin{align}
A_0&:=\int_Me^{-u}\bigg(\big(\overline{\Delta}_gd\Delta_g u,du\big)_g+\frac{5}{12}\Scal_g\(\nabla\Delta_g u,\nabla u\)_g-4\left|\nabla u\right|_g^2\Delta_g^2u\nonumber\\
&\qquad-\frac{2}{3}\Scal_g\left|\nabla u\right|_g^2\Delta_g u-\frac{1}{6}\Scal_g^2\left|\nabla u\right|_g^2\bigg)dv_g\nonumber\\
&\,\,=\(p-4\)\int_Me^{\(p-1\)u}\left|\nabla u\right|_g^2dv_g.\label{Th2Eq2}
\end{align}
Integrating by parts and using \eqref{Sec2Eq1} and \eqref{Sec2Eq2}, we also obtain
\begin{equation}\label{Th2Eq3}
A_1=\dots=A_{12}=0,
\end{equation}
where
\begin{align*}
A_1&:=\int_M e^{-u}\big(\big(\overline{\Delta}_gd\Delta_g u,d u\big)_g-\(\nabla^2\Delta_g u,\nabla^2u-\nabla u\otimes\nabla u\)_g\big)dv_g,\allowdisplaybreaks\\
A_2&:=\int_Me^{-u}\bigg(2\(\nabla^2\Delta_g u,\nabla^2u\)_g-\big(\nabla\Delta_gu,\nabla\left|\nabla u\right|_g^2\big)_g-2\left|\nabla \Delta_g u\right|^2_g\nonumber\\
&\qquad+\frac{1}{2}\Scal_g\(\nabla\Delta_gu,\nabla u\)_g\bigg)dv_g\nonumber\\
&\,\,=\int_Me^{-u}\big(2\(\nabla^2\Delta_g u,\nabla^2u\)_g-\big(\nabla\Delta_gu,\nabla\left|\nabla u\right|_g^2\big)_g-2\big(\overline{\Delta}_g d u,d \Delta_g u\big)_g\big)dv_g,\allowdisplaybreaks\\
A_3&:=\int_M e^{-u}\big(2\(\nabla^2\Delta_g u,\nabla u\otimes\nabla u\)_g+\big(\nabla\Delta_gu,\nabla\left|\nabla u\right|_g^2\big)_g-\(\nabla\Delta_g u,\nabla u\)_g\Delta_g u\nonumber\\
&\qquad-2\left|\nabla u\right|_g^2\(\nabla\Delta_g u,\nabla u\)_g\big)dv_g,\allowdisplaybreaks\\
A_4&:=\int_M e^{-u}\big(\left|\nabla u\right|_g^2\Delta_g^2u-\big(\nabla\Delta_gu,\nabla\left|\nabla u\right|_g^2\big)_g+\left|\nabla u\right|_g^2\(\nabla\Delta_g u,\nabla u\)_g\big)dv_g,\allowdisplaybreaks\\
A_5&:=\int_M e^{-u}\bigg(\big(\nabla\Delta_gu,\nabla\left|\nabla u\right|_g^2\big)_g-2\(\nabla\Delta_g u,\nabla u\)_g\Delta_g u+2\left|\nabla^2 u\right|_g^2\Delta_g u\nonumber\\
&\qquad-\big(\nabla\left|\nabla u\right|_g^2,\nabla u\big)_g\Delta_g u+\frac{1}{2}\Scal_g\left|\nabla u\right|_g^2\Delta_g u\bigg)dv_g\nonumber\\
&\,\,=\int_M e^{-u}\big(\big(\nabla\Delta_gu,\nabla\left|\nabla u\right|_g^2\big)_g-\Delta_g\left|\nabla u\right|_g^2\Delta_g u-\big(\nabla\left|\nabla u\right|_g^2,\nabla u\big)_g\Delta_g u\big)dv_g,\allowdisplaybreaks\\
A_6&:=\int_M e^{-u}\big(2\big(\nabla\Delta_gu,\nabla u\big)_g\Delta_gu-\(\Delta_g u\)^3-\left|\nabla u\right|_g^2\(\Delta_gu\)^2\big)dv_g,\allowdisplaybreaks\\
A_7&:=\int_M e^{-u}\bigg(2\left|\nabla u\right|_g^2\(\nabla\Delta_g u,\nabla u\)_g-2\left|\nabla u\right|_g^2\left|\nabla^2 u\right|_g^2\nonumber\\
&\qquad-\big|\nabla\left|\nabla u\right|_g^2\big|_g^2+\left|\nabla u\right|_g^2\big(\nabla\left|\nabla u\right|_g^2,\nabla u\big)_g-\frac{1}{2}\Scal_g\left|\nabla u\right|_g^4\bigg)dv_g\nonumber\\
&\,\,=\int_M e^{-u}\big(\left|\nabla u\right|_g^2\Delta_g\left|\nabla u\right|_g^2-\big|\nabla\left|\nabla u\right|_g^2\big|_g^2+\left|\nabla u\right|_g^2\big(\nabla\left|\nabla u\right|_g^2,\nabla u\big)_g\big)dv_g,\allowdisplaybreaks\\
A_8&:=\int_M e^{-u}\big(\left|\nabla u\right|_g^2\(\nabla\Delta_g u,\nabla u\)_g+\big(\nabla\left|\nabla u\right|_g^2,\nabla u\big)_g\Delta_g u-\left|\nabla u\right|_g^2\(\Delta_gu\)^2\nonumber\\
&\qquad-\left|\nabla u\right|_g^4\Delta_gu\big)dv_g,\allowdisplaybreaks\\
A_9&:=\int_M e^{-u}\big(2\left|\nabla u\right|_g^2\big(\nabla\left|\nabla u\right|_g^2,\nabla u\big)_g-\left|\nabla u\right|_g^4\Delta_g u-\left|\nabla u\right|_g^6\big)dv_g,\allowdisplaybreaks\\
A_{10}&:=\int_M e^{-u}\bigg(2\(\nabla\Delta_gu,\nabla u\)_g-2\left|\nabla^2u\right|_g^2+\big(\nabla\left|\nabla u\right|_g^2,\nabla u\big)_g-\frac{1}{2}\Scal_g\left|\nabla u\right|_g^2\bigg)dv_g\nonumber\\
&\,\,=\int_M e^{-u}\big(\Delta_g\left|\nabla u\right|_g^2+\big(\nabla\left|\nabla u\right|_g^2,\nabla u\big)_g\big)dv_g,\allowdisplaybreaks\\
A_{11}&:=\int_M e^{-u}\big(\(\nabla\Delta_g u,\nabla u\)_g-\(\Delta_gu\)^2-\left|\nabla u\right|_g^2\Delta_gu\big)dv_g,\allowdisplaybreaks\\
A_{12}&:=\int_M e^{-u}\big(\big(\nabla\left|\nabla u\right|_g^2,\nabla u\big)_g-\left|\nabla u\right|_g^2\Delta_g u-\left|\nabla u\right|_g^4\big)dv_g.
\end{align*}
Combining \eqref{Th2Eq2} and \eqref{Th2Eq3}, we then obtain
\begin{align}\label{Th2Eq15}
&36\(p-4\)\int_Me^{\(p-1\)u}\left|\nabla u\right|_g^2dv_g\nonumber\\
&\quad=36A_0-36A_1-18A_2+18A_3+144A_4+84A_5+42A_6+12A_7-60A_8+18A_9\nonumber\\
&\qquad-20\Scal_gA_{10}+10\Scal_gA_{11}-12\Scal_gA_{12}\allowdisplaybreaks\nonumber\\
&\quad=\int_Me^{-u}\big(36\left|\nabla\Delta_gu\right|^2_g
-24\big(\nabla\Delta_gu,\nabla \left|\nabla u\right|^2_g\big)_g
-120\big(\nabla\Delta_gu,\nabla u\big)_g\Delta_gu\allowdisplaybreaks\nonumber\\
&\qquad+72\left|\nabla u\right|_g^2\(\nabla\Delta_gu,\nabla u\)_g
-12\big|\nabla\left|\nabla u\right|_g^2\big|_g^2
-144\big(\nabla\left|\nabla u\right|_g^2,\nabla u\big)_g\Delta_gu\allowdisplaybreaks\nonumber\\
&\qquad+48\left|\nabla u\right|_g^2\big(\nabla\left|\nabla u\right|_g^2,\nabla u\big)_g
-24\left|\nabla u\right|_g^2\left|\nabla^2 u\right|_g^2
+168\left|\nabla^2 u\right|^2_g\Delta_gu
-42\(\Delta_gu\)^3\allowdisplaybreaks\nonumber\\
&\qquad+18\left|\nabla u\right|_g^2\(\Delta_gu\)^2
+42\left|\nabla u\right|_g^4\Delta_gu
-18\left|\nabla u\right|_g^6
-24\Scal_g\(\nabla\Delta_gu,\nabla u\)_g\allowdisplaybreaks\nonumber\\
&\qquad-32\Scal_g\big(\nabla\left|\nabla u\right|_g^2,\nabla u\big)_g
+40\Scal_g\left|\nabla^2 u\right|^2_g
-10\Scal_g\(\Delta_gu\)^2
+20\Scal_g\left|\nabla u\right|_g^2\Delta_gu\nonumber\\
&\qquad+6\Scal_g\left|\nabla u\right|_g^4
+4\Scal_g^2\left|\nabla u\right|_g^2\big)dv_g.
\end{align}
On the other hand, the transformation law for the scalar curvature under a conformal change of metric gives 
\begin{equation}\label{Th2Eq16}
\Scal_{e^{2u}g}=e^{-3u}\(6\Delta_ge^u+\Scal_ge^u\)=e^{-2u}\big(6\Delta_gu-6\left|\nabla u\right|_g^2+\Scal_g\big).
\end{equation}
Differentiating \eqref{Th2Eq16}, we obtain
\begin{equation}\label{Th2Eq17}
\nabla\Scal_{e^{2u}g}=e^{-2u}\big(6\nabla\Delta_gu-12\Delta_gu\nabla u-6\nabla\left|\nabla u\right|_g^2+12\left|\nabla u\right|_g^2\nabla u-2\Scal_g\nabla u\big).
\end{equation}
Moreover, using \eqref{Sec1Eq3}, we obtain
\begin{align}
\Einstein_{e^{2u}g}&=-2\nabla^2u+2\nabla u\otimes\nabla u-\frac{1}{2}\big(\Delta_g u+\left|\nabla u\right|^2_g\big)g,\label{Th2Eq18}\allowdisplaybreaks\\
\left|\Einstein_{e^{2u}g}\right|^2_g&=4\left|\nabla^2u\right|^2_g-4\big(\nabla\left|\nabla u\right|^2_g,\nabla u\big)_g-\(\Delta_g u\)^2-2\left|\nabla u\right|^2_g\Delta_gu+3\left|\nabla u\right|^4_g\label{Th2Eq19}
\end{align}
and
\begin{multline}\label{Th2Eq20}
\big|\Einstein_{e^{2u}g}\nabla\(e^{-2u}\)\big|_g^2=e^{-4u}\big(
4\big|\nabla\left|\nabla u\right|_g^2\big|_g^2
+4\big(\nabla\left|\nabla u\right|_g^2,\nabla u\big)_g\Delta_gu\\
-12\left|\nabla u\right|_g^2\big(\nabla\left|\nabla u\right|_g^2,\nabla u\big)_g
+\left|\nabla u\right|_g^2\(\Delta_gu\)^2
-6\left|\nabla u\right|_g^4\Delta_gu
+9\left|\nabla u\right|_g^6\big).
\end{multline}
It follows from \eqref{Th2Eq16}--\eqref{Th2Eq20} that
\begin{align}\label{Th2Eq21}
&\big|\nabla\Scal_{e^{2u}g}+2\Einstein_{e^{2u}g}\nabla\(e^{-2u}\)\big|^2_g-4\big|\Einstein_{e^{2u}g}\nabla\(e^{-2u}\)\big|_g^2\nonumber\\
&\quad=4e^{-4u}\big(9\left|\nabla\Delta_gu\right|^2_g
-6\big(\nabla\Delta_gu,\nabla \left|\nabla u\right|^2_g\big)_g
-30\big(\nabla\Delta_gu,\nabla u\big)_g\Delta_gu\nonumber\allowdisplaybreaks\\
&\qquad+18\left|\nabla u\right|_g^2\(\nabla\Delta_gu,\nabla u\)_g
-3\big|\nabla\left|\nabla u\right|_g^2\big|_g^2
+6\big(\nabla\left|\nabla u\right|_g^2,\nabla u\big)_g\Delta_gu\allowdisplaybreaks\nonumber\\
&\qquad+6\left|\nabla u\right|_g^2\big(\nabla\left|\nabla u\right|_g^2,\nabla u\big)_g
+24\left|\nabla u\right|_g^2\(\Delta_gu\)^2
-24\left|\nabla u\right|_g^4\Delta_gu\allowdisplaybreaks\nonumber\\
&\qquad-6\Scal_g\(\nabla\Delta_gu,\nabla u\)_g
+2\Scal_g\big(\nabla\left|\nabla u\right|_g^2,\nabla u\big)_g
+10\Scal_g\left|\nabla u\right|_g^2\Delta_gu\nonumber\\
&\qquad-6\Scal_g\left|\nabla u\right|_g^4
+\Scal_g^2\left|\nabla u\right|_g^2\big)
\end{align}
and
\begin{align}\label{Th2Eq22}
&\left|\Einstein_{e^{2u}g}\right|^2_g\(7\Scal_{e^{2u}g}e^{-2u}+3\Scal_ge^{-4u}+9\left|\nabla e^{-2u}\right|_g^2\)\nonumber\\
&\quad=e^{-4u}\big(-168\big(\nabla\left|\nabla u\right|_g^2,\nabla u\big)_g\Delta_gu
+24\left|\nabla u\right|_g^2\big(\nabla\left|\nabla u\right|_g^2,\nabla u\big)_g\allowdisplaybreaks\nonumber\\
&\qquad-24\left|\nabla u\right|_g^2\left|\nabla^2 u\right|_g^2
+168\left|\nabla^2 u\right|^2_g\Delta_gu
-42\(\Delta_gu\)^3
-78\left|\nabla u\right|_g^2\(\Delta_gu\)^2\allowdisplaybreaks\nonumber\\
&\qquad+138\left|\nabla u\right|_g^4\Delta_gu
-18\left|\nabla u\right|_g^6
-40\Scal_g\big(\nabla\left|\nabla u\right|_g^2,\nabla u\big)_g
+40\Scal_g\left|\nabla^2 u\right|^2_g\nonumber\\
&\qquad-10\Scal_g\(\Delta_gu\)^2
-20\Scal_g\left|\nabla u\right|_g^2\Delta_gu
+30\Scal_g\left|\nabla u\right|_g^4\big).
\end{align}
Putting together \eqref{Th2Eq15}, \eqref{Th2Eq21} and \eqref{Th2Eq22}, we then obtain 
\begin{equation}\label{Th2Eq23}
36\(p-4\)\int_Me^{\(p-1\)u}\left|\nabla u\right|_g^2dv_g=\int_M\Theta^{\(2\)}_g\(e^{-2u}\)dv_g.
\end{equation}
Notice that \eqref{Sec1Eq2} and \eqref{Th2Eq1} give $Q_{e^{2u}g}>0$ in $M$. Moreover, since $\Scal_g>0$ in $M$ and $\(M,g\)$ is Einstein, we obtain that $\Q_g>0$ in $M$. Applying Theorem~\ref{ThGM}, it then follows that $\Scal_{e^{2u}g}>0$ in $M$. Using \eqref{Sec1Eq5} together with the positivity of the functions $v$, $\Scal_g$ and $\Scal_{e^{2u}g}$, we then obtain that $\Theta^{\(2\)}_g\(e^{-2u}\)\ge0$ in $M$. It then follows from \eqref{Th2Eq23} that if $p<4$, then $\nabla u\equiv0$ in $M$ and if $p=4$, then $\Theta^{\(2\)}_g\(e^{-2u}\)\equiv0$ in $M$. Since in the latter case, we assumed that $\(M,g\)$ is not conformally diffeomorphic to the standard sphere, we then obtain that $u$ is constant. This ends the proof of Theorem~\ref{Th2}.
\endproof

We now prove Theorem~\ref{Th3}.

\proof[Proof of Theorem~\ref{Th3}]
Let $u$ be a solution of \eqref{Th2Eq1}. Integrating by parts, we then obtain
\begin{align}
A_0&:=\int_Mu^{-\frac{2}{n-4}}\(\Delta_gu-\frac{n+2}{n-4}u^{-1}\left|\nabla u\right|_g^2\)\Paneitz_gu\,dv_g\nonumber\\
&\,\,=\(p-\frac{2n}{n-4}\)\int_Mu^{p-\frac{2\(n-3\)}{n-4}}\left|\nabla u\right|_g^2dv_g.\label{Th3Eq2}
\end{align}
Integrating by parts and using \eqref{Sec2Eq2} and using the fact that $\(M,g\)$ is Einstein, we also obtain
\begin{equation}\label{Th3Eq3}
A_1=\dots=A_{11}=0,
\end{equation}
where\addtolength{\textheight}{19pt}
\begin{align*}
A_1&:=\int_M u^{-\frac{2}{n-4}}\bigg(\Delta_g^2u\Delta_gu-\left|\nabla\Delta_gu\right|_g^2+\frac{2}{n-4}u^{-1}\(\nabla\Delta_gu,\nabla u\)_g\Delta_gu\bigg)dv_g,\allowdisplaybreaks\\
A_2&:=\int_M u^{-\frac{n-2}{n-4}}\bigg(\left|\nabla u\right|^2_g\Delta_g^2u-\big(\nabla\Delta_gu,\nabla\left|\nabla u\right|^2_g\big)_g\nonumber\\
&\qquad+\frac{n-2}{n-4}u^{-1}\left|\nabla u\right|^2_g\(\nabla\Delta_gu,\nabla u\)_g\bigg)dv_g,\allowdisplaybreaks\\
A_3&:=\int_M u^{-\frac{n-2}{n-4}}\bigg(\big(\nabla\Delta_gu,\nabla\left|\nabla u\right|^2_g\big)_g-2\(\nabla\Delta_g u,\nabla u\)_g\Delta_gu+2\left|\nabla^2 u\right|_g^2\Delta_gu\nonumber\\
&\qquad-\frac{n-2}{n-4}u^{-1}\big(\nabla\left|\nabla u\right|^2_g,\nabla u\big)_g\Delta_gu+\frac{2}{n}\Scal_g\left|\nabla u\right|^2_g\Delta_gu\bigg)dv_g\nonumber\\
&\,\,=\int_M u^{-\frac{n-2}{n-4}}\bigg(\big(\nabla\Delta_gu,\nabla\left|\nabla u\right|^2_g\big)_g-\Delta_g\left|\nabla u\right|_g^2\Delta_gu\nonumber\\
&\qquad-\frac{n-2}{n-4}u^{-1}\big(\nabla\left|\nabla u\right|^2_g,\nabla u\big)_g\Delta_gu\bigg)dv_g,\allowdisplaybreaks\\
A_4&:=\int_M u^{-\frac{n-2}{n-4}}\bigg(2\(\nabla\Delta_g u,\nabla u\)_g\Delta_gu-\(\Delta_gu\)^3-\frac{n-2}{n-4}u^{-1}\left|\nabla u\right|^2_g\(\Delta_gu\)^2\bigg)dv_g,\allowdisplaybreaks\\
A_5&:=\int_M u^{-\frac{2\(n-3\)}{n-4}}\bigg(2\left|\nabla u\right|^2_g\(\nabla\Delta_gu,\nabla u\)_g-2\left|\nabla u\right|^2_g\left|\nabla^2 u\right|^2_g-\big|\nabla\left|\nabla u\right|^2_g\big|^2_g\nonumber\\
&\qquad+\frac{2\(n-3\)}{n-4}u^{-1}\left|\nabla u\right|^2_g\big(\nabla\left|\nabla u\right|^2_g,\nabla u\big)_g-\frac{2}{n}\Scal_g\left|\nabla u\right|^4_g\bigg)dv_g\nonumber\\
&\,\,=\int_M u^{-\frac{2\(n-3\)}{n-4}}\bigg(\left|\nabla u\right|^2_g\Delta_g\left|\nabla u\right|^2_g-\big|\nabla\left|\nabla u\right|^2_g\big|^2_g\nonumber\\
&\qquad+\frac{2\(n-3\)}{n-4}u^{-1}\left|\nabla u\right|^2_g\big(\nabla\left|\nabla u\right|^2_g,\nabla u\big)_g\bigg)dv_g,\allowdisplaybreaks\\
A_6&:=\int_M u^{-\frac{2\(n-3\)}{n-4}}\bigg(\left|\nabla u\right|^2_g\(\nabla\Delta_g u,\nabla u\)_g+\big(\nabla\left|\nabla u\right|^2_g,\nabla u\big)_g\Delta_g u-\left|\nabla u\right|_g^2\(\Delta_gu\)^2\nonumber\\
&\qquad-\frac{2\(n-3\)}{n-4}u^{-1}\left|\nabla u\right|^4_g\Delta_g u\bigg)dv_g,\allowdisplaybreaks\\
A_7&:=\int_M u^{-\frac{3n-10}{n-4}}\bigg(2\left|\nabla u\right|^2_g\big(\nabla\left|\nabla u\right|^2_g,\nabla u\big)_g-\left|\nabla u\right|^4_g\Delta_gu\nonumber\\
&\qquad-\frac{3n-10}{n-4}u^{-1}\left|\nabla u\right|^6_g\bigg)dv_g,\allowdisplaybreaks\\
A_8&:=\int_M u^{-\frac{2}{n-4}}\bigg(\(\nabla\Delta_gu,\nabla u\)_g-\(\Delta_gu\)^2-\frac{2}{n-4}u^{-1}\left|\nabla u\right|_g^2\Delta_gu\bigg)dv_g,\allowdisplaybreaks\\
A_9&:=\int_M u^{-\frac{2}{n-4}}\bigg(2\(\nabla\Delta_g u,\nabla u\)_g-2\left|\nabla^2 u\right|_g^2+\frac{2}{n-4}u^{-1}\big(\nabla\left|\nabla u\right|^2_g,\nabla u\big)_g\nonumber\\
&\quad-\frac{2}{n}\Scal_g\left|\nabla u\right|^2_g\bigg)dv_g\nonumber\\
&\,\,=\int_M u^{-\frac{2}{n-4}}\bigg(\Delta_g\left|\nabla u\right|_g^2+\frac{2}{n-4}u^{-1}\big(\nabla\left|\nabla u\right|^2_g,\nabla u\big)_g\bigg)dv_g,\allowdisplaybreaks\\
A_{10}&:=\int_Mu^{-\frac{n-2}{n-4}}\bigg(\big(\nabla\left|\nabla u\right|^2_g,\nabla u\big)_g-\left|\nabla u\right|^2_g\Delta_gu-\frac{n-2}{n-4}u^{-1}\left|\nabla u\right|^4_g\bigg)dv_g,\allowdisplaybreaks\\
A_{11}&:=\int_M u^{\frac{n-6}{n-4}}\bigg(\Delta_gu-\frac{n-6}{n-4}u^{-1}\left|\nabla u\right|^2_g\bigg)dv_g.
\end{align*}
Combining \eqref{Th3Eq2} and \eqref{Th3Eq3}, we then obtain\addtolength{\textheight}{18pt}
\begin{align}\label{Th3Eq14}
&\frac{16\(n-1\)^2}{\(n-4\)^2}\(p-\frac{2n}{n-4}\)\int_Mu^{p-\frac{2\(n-3\)}{n-4}}\left|\nabla u\right|_g^2dv_g\nonumber\\
&\quad=\frac{16\(n-1\)^2}{\(n-4\)^2}A_0-\frac{16\(n-1\)^2}{\(n-4\)^2}A_1+\frac{16\(n-1\)^2\(n+2\)}{\(n-4\)^3}A_2\allowdisplaybreaks\nonumber\\
&\qquad+\frac{16\(n-1\)\(n^2-2\)}{\(n-4\)^3}A_3+\frac{32\(n-1\)\(n^2-2\)}{n\(n-4\)^3}A_4+\frac{32\(n-1\)\(n-2\)}{\(n-4\)^4}A_5\allowdisplaybreaks\nonumber\\
&\qquad-\frac{16\(n-1\)\(n-2\)\(n^3-n^2-4n+8\)}{n\(n-4\)^4}A_6+\frac{64\(n-1\)^2\(n-2\)^2}{n\(n-4\)^5}A_7\allowdisplaybreaks\nonumber\\
&\qquad+\frac{8n\(n-2\)}{\(n-4\)^2}\Scal_gA_8-\frac{4\(n^2+2n-4\)}{\(n-4\)^2}\Scal_gA_9-\frac{8\(n-1\)\(n^2-12\)}{\(n-4\)^3}\Scal_gA_{10}\nonumber\\
&\qquad-\frac{\(n-2\)\(n+2\)}{n\(n-4\)}\Scal_g^2A_{11}\nonumber\allowdisplaybreaks\\
&\quad=\frac{8}{\(n-4\)^2}\int_Mu^{-\frac{2}{n-4}}\bigg(2\(n-1\)^2\left|\nabla\Delta_gu\right|^2_g
-\frac{2n\(n-1\)}{n-4}u^{-1}\big(\nabla\Delta_gu,\nabla \left|\nabla u\right|^2_g\big)_g\allowdisplaybreaks\nonumber\\
&\qquad-\frac{4\(n-1\)\(n^3-n^2-3n+4\)}{n\(n-4\)}u^{-1}\big(\nabla\Delta_gu,\nabla u\big)_g\Delta_gu
\allowdisplaybreaks\nonumber\\
&\qquad+\frac{4\(n-1\)^2\(n-2\)\(n+4\)}{n\(n-4\)^2}u^{-2}\left|\nabla u\right|_g^2\(\nabla\Delta_gu,\nabla u\)_g\allowdisplaybreaks\nonumber\\
&\qquad-\frac{4\(n-1\)\(n-2\)}{\(n-4\)^2}u^{-2}\big|\nabla\left|\nabla u\right|_g^2\big|_g^2
\allowdisplaybreaks\nonumber\\
&\qquad-\frac{2\(n-1\)\(n-2\)\(n+2\)\(2n^2-5n+4\)}{n\(n-4\)^2}u^{-2}\big(\nabla\left|\nabla u\right|_g^2,\nabla u\big)_g\Delta_gu\allowdisplaybreaks\nonumber\\
&\qquad+\frac{8\(n-1\)\(n-2\)\(3n^2-9n+4\)}{n\(n-4\)^3}u^{-3}\left|\nabla u\right|_g^2\big(\nabla\left|\nabla u\right|_g^2,\nabla u\big)_g\allowdisplaybreaks\nonumber\\
&\qquad-\frac{8\(n-1\)\(n-2\)}{\(n-4\)^2}u^{-2}\left|\nabla u\right|_g^2\left|\nabla^2 u\right|_g^2
+\frac{4\(n-1\)\(n^2-2\)}{n-4}u^{-1}\left|\nabla^2 u\right|^2_g\Delta_gu
\allowdisplaybreaks\nonumber\\
&\qquad-\frac{4\(n-1\)\(n^2-2\)}{n\(n-4\)}u^{-1}\(\Delta_gu\)^3
\allowdisplaybreaks\nonumber\\
&\qquad+\frac{2\(n-1\)\(n-2\)^2\(n-3\)\(n+2\)}{n\(n-4\)^2}u^{-2}\left|\nabla u\right|_g^2\(\Delta_gu\)^2\allowdisplaybreaks\nonumber\\
&\qquad+\frac{4\(n-1\)\(n-2\)\(n^4-4n^3-3n^2+26n-28\)}{n\(n-4\)^3}u^{-3}\left|\nabla u\right|_g^4\Delta_gu\allowdisplaybreaks\nonumber\\
&\qquad-\frac{8\(n-1\)^2\(n-2\)^2\(3n-10\)}{n\(n-4\)^4}u^{-4}\left|\nabla u\right|_g^6
-4\(n-1\)\Scal_g\(\nabla\Delta_gu,\nabla u\)_g\allowdisplaybreaks\nonumber\\
&\qquad-\frac{n^3-10n+8}{n-4}\Scal_gu^{-1}\big(\nabla\left|\nabla u\right|_g^2,\nabla u\big)_g
+\(n^2+2n-4\)\Scal_g\left|\nabla^2 u\right|^2_g\allowdisplaybreaks\nonumber\\
&\qquad-\frac{n^2+2n-4}{n}\Scal_g\(\Delta_gu\)^2
+\frac{2\(n^2-2n+2\)}{n-4}\Scal_gu^{-1}\left|\nabla u\right|_g^2\Delta_gu\nonumber\\
&\qquad+\frac{\(n-1\)\(n-2\)\(n^3-12n-8\)}{n\(n-4\)^2}\Scal_gu^{-2}\left|\nabla u\right|_g^4
+2\Scal_g^2\left|\nabla u\right|_g^2\bigg)dv_g.
\end{align}
On the other hand, the transformation law for the scalar curvature under a conformal change of metric gives\addtolength{\textheight}{-37pt}
\begin{align}\label{Th3Eq15}
\Scal_{u^{\frac{4}{n-4}}g}&=u^{-\frac{n+2}{n-4}}\bigg(\frac{4\(n-1\)}{n-2}\Delta_g\(u^{\frac{n-2}{n-4}}\)+\Scal_gu^{\frac{n-2}{n-4}}\bigg)\nonumber\\
&=u^{-\frac{n}{n-4}}\bigg(\frac{4\(n-1\)}{n-4}\Delta_gu-\frac{8\(n-1\)}{\(n-4\)^2}u^{-1}\left|\nabla u\right|_g^2+\Scal_gu\bigg).
\end{align}
Differentiating \eqref{Th3Eq15}, we obtain
\begin{multline}\label{Th3Eq16}
\nabla\Scal_{u^{\frac{4}{n-4}}g}=\frac{4}{n-4}u^{-\frac{n}{n-4}}\bigg(\(n-1\)\nabla\Delta_gu-\frac{n\(n-1\)}{n-4}u^{-1}\Delta_gu\nabla u\\
-\frac{2\(n-1\)}{n-4}u^{-1}\nabla\left|\nabla u\right|_g^2+\frac{4\(n-1\)\(n-2\)}{\(n-4\)^2}u^{-2}\left|\nabla u\right|_g^2\nabla u-\Scal_g\nabla u\bigg).
\end{multline}
Moreover, using \eqref{Sec1Eq3}, we obtain
\begin{align}
\Einstein_{u^{\frac{4}{n-4}}g}&=-\frac{2\(n-2\)}{n-4}u^{-1}\bigg(\nabla^2u-\frac{n-2}{n-4}u^{-1}\nabla u\otimes\nabla u\nonumber\\
&\qquad+\frac{1}{n}\(\Delta_gu+\frac{n-2}{n-4}u^{-1}\left|\nabla u\right|^2_g\)g\bigg),\label{Th3Eq17}\allowdisplaybreaks\\
\left|\Einstein_{u^{\frac{4}{n-4}}g}\right|^2_g&=\frac{4\(n-2\)^2}{\(n-4\)^2}u^{-2}\bigg(\left|\nabla^2u\right|^2_g-\frac{n-2}{n-4}u^{-1}\big(\nabla\left|\nabla u\right|^2_g,\nabla u\big)_g-\frac{1}{n}\(\Delta_g u\)^2\nonumber\\
&\qquad-\frac{2\(n-2\)}{n\(n-4\)}u^{-1}\left|\nabla u\right|^2_g\Delta_gu+\frac{\(n-1\)\(n-2\)^2}{n\(n-4\)^2}u^{-2}\left|\nabla u\right|^4_g\bigg)\label{Th3Eq18}
\end{align}
and
\begin{align}\label{Th3Eq19}
&\left|\Einstein_{u^{\frac{4}{n-4}}g}\nabla\(u^{-\frac{4}{n-4}}\)\right|_g^2=\frac{16\(n-2\)^2}{\(n-4\)^4}u^{-\frac{4\(n-2\)}{n-4}}\bigg(
\big|\nabla\left|\nabla u\right|_g^2\big|_g^2\nonumber\\
&\qquad+\frac{4}{n}\big(\nabla\left|\nabla u\right|_g^2,\nabla u\big)_g\Delta_gu-\frac{4\(n-1\)\(n-2\)}{n\(n-4\)}u^{-1}\left|\nabla u\right|_g^2\big(\nabla\left|\nabla u\right|_g^2,\nabla u\big)_g
\allowdisplaybreaks\nonumber\\
&\qquad+\frac{4}{n^2}\left|\nabla u\right|_g^2\(\Delta_gu\)^2-\frac{8\(n-1\)\(n-2\)}{n^2\(n-4\)}u^{-1}\left|\nabla u\right|_g^4\Delta_gu
\nonumber\\
&\qquad+\frac{4\(n-1\)^2\(n-2\)^2}{n^2\(n-4\)^2}u^{-2}\left|\nabla u\right|_g^6\bigg).
\end{align}
It follows from \eqref{Th3Eq15}--\eqref{Th3Eq19} that
\begin{align}\label{Th3Eq20}
&\left|\nabla\Scal_{u^{\frac{4}{n-4}}g}+\frac{3n-4}{2\(n-2\)}\Einstein_{u^{\frac{4}{n-4}}g}\nabla\(u^{-\frac{4}{n-4}}\)\right|^2_g-\frac{\(3n-4\)^2}{4\(n-2\)^2}\left|\Einstein_{u^{\frac{4}{n-4}}g}\nabla\(u^{-\frac{4}{n-4}}\)\right|_g^2\nonumber\\
&\quad=\frac{16}{\(n-4\)^2}u^{-\frac{2n}{n-4}}\bigg(\(n-1\)^2\left|\nabla\Delta_gu\right|^2_g
-\frac{n\(n-1\)}{n-4}u^{-1}\big(\nabla\Delta_gu,\nabla \left|\nabla u\right|^2_g\big)_g\allowdisplaybreaks\nonumber\\
&\qquad-\frac{2\(n-1\)\(n^3-n^2-3n+4\)}{n\(n-4\)}u^{-1}\big(\nabla\Delta_gu,\nabla u\big)_g\Delta_gu
\allowdisplaybreaks\nonumber\\
&\qquad+\frac{2\(n-1\)^2\(n-2\)\(n+4\)}{n\(n-4\)^2}u^{-2}\left|\nabla u\right|_g^2\(\nabla\Delta_gu,\nabla u\)_g\allowdisplaybreaks\nonumber\\
&\qquad-\frac{2\(n-1\)\(n-2\)}{\(n-4\)^2}u^{-2}\big|\nabla\left|\nabla u\right|_g^2\big|_g^2\allowdisplaybreaks\nonumber\\
&\qquad+\frac{\(n-1\)\(n-2\)^2\(n+4\)}{n\(n-4\)^2}u^{-2}\big(\nabla\left|\nabla u\right|_g^2,\nabla u\big)_g\Delta_gu\allowdisplaybreaks\nonumber\\
&\qquad+\frac{4\(n-1\)\(n-2\)\(2n^2-7n+4\)}{n\(n-4\)^3}u^{-3}\left|\nabla u\right|_g^2\big(\nabla\left|\nabla u\right|_g^2,\nabla u\big)_g\allowdisplaybreaks\nonumber\\
&\qquad+\frac{\(n-1\)\(n-2\)\(n^2+n-4\)}{\(n-4\)^2}u^{-2}\left|\nabla u\right|_g^2\(\Delta_gu\)^2\allowdisplaybreaks\nonumber\\
&\qquad-\frac{2\(n-1\)\(n-2\)\(n^3+3n^2-16n+16\)}{n\(n-4\)^3}u^{-3}\left|\nabla u\right|_g^4\Delta_gu\allowdisplaybreaks\nonumber\\
&\qquad-\frac{8\(n-1\)^2\(n-2\)^2}{n\(n-4\)^3}u^{-4}\left|\nabla u\right|_g^6
-2\(n-1\)\Scal_g\(\nabla\Delta_gu,\nabla u\)_g\allowdisplaybreaks\nonumber\\
&\qquad+\frac{n}{n-4}\Scal_gu^{-1}\big(\nabla\left|\nabla u\right|_g^2,\nabla u\big)_g
+\frac{2\(n^3-n^2-3n+4\)}{n\(n-4\)}\Scal_gu^{-1}\left|\nabla u\right|_g^2\Delta_gu\nonumber\\
&\qquad-\frac{2\(n-1\)\(n-2\)\(n+4\)}{n\(n-4\)^2}\Scal_gu^{-2}\left|\nabla u\right|_g^4
+\Scal_g^2\left|\nabla u\right|_g^2\bigg)
\end{align}
and
\begin{align}\label{Th3Eq21}
&\left|\Einstein_{u^{\frac{4}{n-4}}g}\right|^2_g\big(2\(n^2-2\)\Scal_{e^{2u}g}u^{-\frac{4}{n-4}}+4\(n-1\)\Scal_gu^{-\frac{8}{n-4}}+n\(n-1\)^2\big|\nabla\big(u^{-\frac{4}{n-4}}\big)\big|_g^2\big)\nonumber\\
&\quad=\frac{8\(n-2\)^2}{\(n-4\)^2}u^{-\frac{2n}{n-4}}\bigg(-\frac{4\(n-1\)\(n-2\)\(n^2-2\)}{\(n-4\)^2}u^{-2}\big(\nabla\left|\nabla u\right|_g^2,\nabla u\big)_g\Delta_gu\allowdisplaybreaks\nonumber\\
&\qquad+\frac{8\(n-1\)\(n-2\)^2}{\(n-4\)^3}u^{-3}\left|\nabla u\right|_g^2\big(\nabla\left|\nabla u\right|_g^2,\nabla u\big)_g\allowdisplaybreaks\nonumber\\
&\qquad-\frac{8\(n-1\)\(n-2\)}{\(n-4\)^2}u^{-2}\left|\nabla u\right|_g^2\left|\nabla^2 u\right|_g^2\allowdisplaybreaks\nonumber\\
&\qquad+\frac{4\(n-1\)\(n^2-2\)}{n-4}u^{-1}\left|\nabla^2 u\right|^2_g\Delta_gu
-\frac{4\(n-1\)\(n^2-2\)}{n\(n-4\)}u^{-1}\(\Delta_gu\)^3\allowdisplaybreaks\nonumber\\
&\qquad-\frac{8\(n-1\)\(n-2\)\(n^2-3\)}{n\(n-4\)^2}u^{-2}\left|\nabla u\right|_g^2\(\Delta_gu\)^2\allowdisplaybreaks\nonumber\\
&\qquad+\frac{4\(n-1\)\(n-2\)^2\(n^3-n^2-2n+6\)}{n\(n-4\)^3}u^{-3}\left|\nabla u\right|_g^4\Delta_gu\allowdisplaybreaks\nonumber\\
&\qquad-\frac{8\(n-1\)^2\(n-2\)^3}{n\(n-4\)^4}u^{-4}\left|\nabla u\right|_g^6\allowdisplaybreaks\nonumber\\
&\qquad-\frac{\(n-2\)\(n^2+2n-4\)}{n-4}\Scal_gu^{-1}\big(\nabla\left|\nabla u\right|_g^2,\nabla u\big)_g\allowdisplaybreaks\nonumber\\
&\qquad+\(n^2+2n-4\)\Scal_g\left|\nabla^2 u\right|^2_g
-\frac{n^2+2n-4}{n}\Scal_g\(\Delta_gu\)^2\allowdisplaybreaks\nonumber\\
&\qquad-\frac{2\(n-2\)\(n^2+2n-4\)}{n\(n-4\)}\Scal_gu^{-1}\left|\nabla u\right|_g^2\Delta_gu\nonumber\\
&\qquad+\frac{\(n-1\)\(n-2\)^2\(n^2+2n-4\)}{n\(n-4\)^2}\Scal_gu^{-2}\left|\nabla u\right|_g^4\bigg)
\end{align}
Putting together \eqref{Th3Eq14}, \eqref{Th3Eq20} and \eqref{Th3Eq21}, we then obtain
$$\frac{16\(n-1\)^2}{\(n-4\)^2}\(p-\frac{2n}{n-4}\)\int_Mu^{p-\frac{2\(n-3\)}{n-4}}\left|\nabla u\right|_g^2dv_g=\int_M\Theta^{\(2\)}_g\(u^{-\frac{4}{n-4}}\)dv_g.$$
We then conclude as in the proof of Theorem~\ref{Th2}.
\endproof

\end{document}